\documentclass[twocolumn,showpacs,preprintnumbers,amsmath,amssymb,showkeys]{revtex4}

\usepackage{graphicx}
\usepackage{dcolumn}
\usepackage{bm}

\newcommand{\eps}{\varepsilon}
\newcommand{\ph}{\varphi}

\newcommand{\trn}{^{\rm\scriptscriptstyle T}}

\newcommand{\mR}{\mathbb R}

\newcommand{\inr}{\!\in \mR}

\DeclareMathOperator{\sign}{sign}

\newcommand{\be}{\begin{equation}}
\newcommand{\ee}{\end{equation}}

\renewcommand{\dfrac}{\frac}

\begin{document}


\title{Controlled Synchronization of One Class of Nonlinear Systems under Information Constraints}
\author{Alexander L. Fradkov, Boris Andrievsky}
  \email{fradkov@mail.ru, bandri@yandex.ru}%
\affiliation{Institute for Problems of Mechanical Engineering,
Russian Academy of Sciences, \\ 61, Bolshoy V.O. Av., 199178,
Saint Petersburg, Russia %
}%
\author{Robin J. Evans}
\email{r.evans@ee.unimelb.edu.au.}%
\affiliation{National ICT Australia,\\%
Department of Electrical and electronic Engineering,\\%
University of Melbourne, Victoria, 3010, Australia}%
\date{\today}

\begin{abstract}
Output feedback controlled synchronization problems for a class of nonlinear unstable systems under
information constraints imposed by limited capacity of
the communication channel are analyzed. 
A binary time-varying coder-decoder scheme is described and
 a theoretical analysis for 
multi-dimensional master-slave systems represented in Lurie
form (linear part plus nonlinearity depending only on measurable
outputs) is provided. An output feedback control law is proposed based on the Passification Theorem.
It is shown that the synchronization error exponentially tends to zero for  sufficiantly high
 transmission rate (channel capacity).  The results obtained for synchronization problem can be extended to tracking problems in a straightforward manner, if the reference signal is described by an {external} ({exogenious}) state space model.
 
 The results are
applied to controlled synchronization of two chaotic Chua systems
 via a communication channel with limited capacity.
 \end{abstract}

\pacs{05.45.Xt, 05.45.Gg}
\keywords{Nonlinear systems, Synchronization, Control,  Communication constraints}

\maketitle
\section{Introduction}

Analysis and control of the behavior of complex interconnected systems and networks has attracted considerable recent interest. The available results significantly depend on models of interconnection between nodes. In some works the interconnections are modeled as delay elements. However, the spatial separation between nodes means that modeling connections via communication channels with limited capacity is more realistic. 

Recently the limitations of control under constraints imposed by a finite capacity information channel have been investigated in detail in the control literature, see \cite{WongBrockett_AC97e,NairEvans_CDC02,NairEvans_Aut03,%
NairEvans_SIAM04,NairEvans_AC04,BazziMitter_InTh05,NairFagnani07} and the references therein. It has been shown that stabilization of linear systems under information constraints is possible if and only if the capacity of the information channel exceeds the entropy production of the system at the equilibrium ({\it Data Rate Theorem}) \cite{NairEvans_Aut03,NairEvans_SIAM04,NairEvans_AC04}. In \cite{Lloyd00,Lloyd04} a general statement was proposed, claiming that the difference between the entropies of the open loop and the closed loop systems cannot exceed the information introduced by the controller, including the transmission rate of the information channel.

For nonlinear systems only a few results are available in the literature \cite{liberzon03,NairEvans_AC04,DePersis_AC05,DePersis_IJRNC06,DePersisNesic_CDCECC05,SavkinCheng_AC07}. In the above papers only the problems of stabilization to a point are considered.

In the control literature there is a strong interest in control of oscillations, particularly in controlled synchronization problems \cite{SpecIssue97,LuChen05,ChopraSpong06,LoriaZavala07}.  However, results of the previous works on control systems analysis under information constraints do not apply to synchronization systems since in a synchronization problem trajectories in the phase space converge to a set (a manifold) rather than to a point, i.e. the problem cannot be reduced to simple stabilization. Moreover, the Data Rate Theorem is difficult to extend to nonlinear systems.

The first results on synchronization under information constraints were presented in \cite{FradkovAndrievskyEvans_PRE06},
where the so called observer-based synchronization scheme \cite{PecoraCarroll90,FNM00} was considered. In this paper we extend the results of  \cite{FradkovAndrievskyEvans_PRE06} and analyze an output feedback controlled synchronization scheme for two nonlinear systems. A major difficulty with the controlled synchronization problem arises because the coupling is implemented in a restricted manner via the control signal which is computed based on a measurable innovation (error) signal which has been transmitted over a communication channel.
 Key tools used to solve the problem are quadratic Lyapunov functions and passification methods \cite{Fradkov74,FradkovMN99}. 
To minimize technicalities we restrict our analysis to Lurie systems (linear part plus nonlinearity depending only on measurable
outputs).

The paper is organized as follows. The controlled synchronization problem is described in Section~\ref{Sec:1}. The coding procedure used in the paper, is presented in Sec.~\ref{Sec:coding}. The main results are presented in Section~\ref{Sec:eval} where an exponential convergence of the synchronization error to zero is established. An example showing synchronization of the chaotic Chua systems is presented in Section \ref{Sec:example}. Final remarks are given in the Conclusion. Auxiliary results are presented in the Appendices.

 \section{Description of controlled synchronization scheme}\label{Sec:1}

Consider two identical dynamical systems modeled in Lurie form (i.e. the right hand sides are split into a linear part and a nonlinear part which depends only on the measurable outputs). Let one of the systems be controlled by a scalar control function $u(t)$ whose action is restricted by a vector of control efficiencies $B$. The controlled system model is as follows:
\begin{align} 
&\dot x(t)=Ax(t)+B\ph (y_1),~y_1(t)=Cx(t), \label{2}
\\
&\dot z(t)=Az(t)+B\ph (y_2)+Bu,~y_2(t)=Cz(t), \label{1}
\end{align} 
where $x(t)$, $z(t)$ are $n$-dimensional (column) vectors of state
variables; $y_1(t)$, $y_2(t)$ are scalar output variables; $A$ is 
an $(n\times n)$-matrix; $B$ is $n\times 1$ (column) matrix; $C$ is  an $1\times n$ (row) matrix, $\ph(y)$ is a continuous nonlinearity, acting in the span of control; vectors $\dot x$, $\dot z $ stand for time-derivatives of $x(t)$, $z(t)$ respectively. System (\ref{2}) is called the {\it master} ({\it leader}) {\it system}, while the controlled system (\ref{1}) is called the {\it slave (follower) system}. Our goal is to evaluate limitations imposed on the synchronization precision by limiting the transmission rate between the systems. The intermediate problem is to find a control function ${\cal U}(\cdot)$ depending on the measurable variables such that the  synchronization error $e(t),$ where $e(t)=x(t)-z(t)$ becomes small as $t$ becomes large. We are also interested in the value of the output synchronization error $\eps(t)=y_1(t) -y_2(t)=Ce(t)$.

A key difficulty arises because the error signal between the master system and the slave systems is not available directly but only through a communication channel with a limited capacity. This means that the synchronization error $\eps(t)$ must be coded at the transmitter side and codewords then transmitted with only a finite number of symbols per second thus introducing error. We assume that the observed signal $\eps(t))$ is coded with symbols from a finite alphabet at discrete sampling time instants $t_k=k T_s$, $k=0,1,2,\dots$, where $T_s$ is the sampling time. Let the coded symbol $\bar \eps [k]=\bar \eps (t_k)$ be
transmitted over a digital communication channel with a finite capacity. To simplify the analysis,
we assume that the observations are not corrupted by observation noise; transmission delay 
and transmission channel distortions may be neglected. Therefore, the discrete communication channel with  sampling period $T_s$ is considered, but it is assumed that the coded symbols are available at the receiver side at the same sampling instant $t_k=kT_s$, as they are generated by the coder.
Assume that {\it zero-order extrapolation} is used to convert the digital sequence $\bar\eps [k]$ to the continuous-time input of the controller $\bar \eps(t)$, namely, that $\bar \eps(t)=\bar\eps[k]$ as $ kT_s \le t< (k+1)T_s $. Then the {\it transmission error} is defined as follows: 
\begin{align}
\delta_\eps(t) = \eps (t)- \bar \eps (t).
\label{dy}
\end{align}
On the receiver side the signal is decoded introducing additional error and the controller can use only the signal $\bar \eps (t)= \eps (t)-\delta_\eps (t)$ instead of $\eps (t)$. 

We restrict consideration to simple control functions in the form of static linear feedback
\begin{align} 
u(t)= K\eps(t), \label{3}
\end{align} 
where $\eps(t)=y_1(t)-y_2(t)$ denotes an output synchronization error and $K$ is a scalar controller gain. The problem of finding static output feedback 
even for linear systems is one of the classical problems of control theory. Although substantial effort has been devoted to its solution and various necessary and sufficient conditions for stabilizability by static output feedback have been obtained, most existing conditions are not testable practically \cite{Syrmos97,Stefanovski06}. In this paper we analyze a natural and relatively broad class of systems for which constructive conditions for output feedback stabilization are known is the class of {\it passifiable} (or {\it feedback passive}) systems (for linear systems this introduced and studied in \cite{Fradkov74,FradkovEJC}). Since we are dealing with a nonlinear problem further complicated by information constraints, we restrict our attention to sufficient conditions for solvability of the problem and evaluate upper bounds for synchronization error. 

\section{Coding procedures}\label{Sec:coding}

In the paper \cite{FradkovAndrievskyEvans_PRE06} the properties of observer-based synchronization for Lurie systems over a limited data rate communication channel with a one-step memory time-varying coder are studied.  It is shown that an upper bound on the limit synchronization error is proportional to a certain upper bound on the transmission error. Under the assumption that a sampling time may be properly chosen, optimality of binary coding in the sense of demanded transmission rate is established, and the relationship between synchronization accuracy and an optimal sampling time is found. On the basis of these results, the present paper deals with a binary coding procedure.

Consider the memoryless (static) binary quantizer to be a discretized map $q: \mR\to\mR$ as 
\begin{align}
\label{qu1}
q (y,M)=M\sign(y), 
\end{align}
where $\sign(\cdot)$ is the {\it signum} function: $\sign(y)=1$, if $y\ge 0$, $\sign (y)=-1$, if $y<0$. Parameter $M$ may be referred to as the {\it quantizer range}. Notice that for a binary coder each codeword symbol contains one bit of information. The discretized output of the considered quantizer is given as $\bar y=q(y,M)$. We assume that the coder and decoder make decisions based on the same information. The output signal of the quantizer is represented as a one-bit information symbol from the coding alphabet ${\cal S}$ and transmitted over the communication channel to the decoder.

In {\it time-varying quantizers}  \cite{BrockettLiberzon_AC00,%
liberzon03,TatikondaMitter_AC04,FradkovAndrievskyEvans_PRE06,NairFagnani07} the range $M$ is updated with time and different values of $M$ are used at each step, $M=M[k]$. Using such a ``zooming'' strategy it is possible to increase coder accuracy in the steady-state mode and at the same time, to prevent coder saturation at the beginning of the process \cite{BrockettLiberzon_AC00}.

In the present paper we use the following time-based zooming strategy for a quantizer range
\begin{align}
\label{mk}
M[k]= M_0\rho^k, ~~k=0,1,\dots ,
\end{align}
where $0<\rho\le 1$ is the decay parameter. The initial value $M_0$ should be large enough to capture the region of possible initial values of $y_0$. Equations \eqref{qu1},  \eqref{mk} describe the coder algorithm. A similar algorithm is realized by the decoder. Namely, the sequence $M[k]$ is reproduced at the receiver node utilizing \eqref{mk} such that the values of $\bar y[k]$ are restored with the given $M[k]$ using the received codeword $s[k]\in{\cal S}$.

\section{Evaluation of synchronization error}\label{Sec:eval}
Let us evaluate the limit synchronization error, taking into account transmission of the error signal over the communication channel and coding procedure. Since the control signal is piecewise constant over sampling intervals $[t_k,t_{k+1}]$, the control law \eqref{3} becomes
\begin{align} 
u(t)= K\bar{\eps}(t), \label{3bar}
\end{align} 
where  $ \bar{\eps}(t)=\bar{\eps}[k]$ as $ k T_s<t<(k+1)T_s$, $\bar{\eps}[k]$ is the result of transmission of the synchronization error signal $\eps(t)=y_1(t)-y_2(t)$ over the channel, $k=0,1,\dots$. 

According to the quantization algorithm \eqref{qu1}, the quantized error signal $\bar{\eps}[k]$ becomes
\begin{align}
\label{qu11}
\bar{\eps}[k]=M[k]\sign(\eps(t_k)), 
\end{align}
where $\eps(t)=y_1(t)-y_2(t)$, $t_k=kT_s$, and the range $M[k]$ is defined by \eqref{mk}. 

Taking into account the stepwise shape of the control function in \eqref{3bar}, rewrite the controller model in the following form:
\begin{align}
u(t)= K\eps(t)-K\delta (t),
\label{02}
\end{align}
where $ \delta (t)=\delta_q (t)+\delta_s(t)$ is a {\it total error}, $\delta_q(t)= \eps(t_k)-\bar{\eps}[k] = Ce(t_k)-\bar{\eps}[k] $ is a {\it quantization error}, $\delta_s(t)=\eps(t)-\eps(t_k)=Ce(t)-Ce(t_k)$ is a {\it sampling error}.

It is seen from the quantization procedure \eqref{qu11} that if the value $\eps(t_k)$ satisfies the inequality $|\eps(t_k)|\le 2M[k]$, then the quantization error does not exceed $M[k]$: $\delta_q(t)|\le M [k]$. For the sampling error the following bounds hold:
\begin{align*}
|\delta_s(t)|\le \int\limits_{t_k}^t{|\dot{\eps}(\tau)|\,{\rm d} \tau}\le\max\limits_{t_k\le\tau\le t}|\dot{\eps}(\tau)|\cdot T_s.
\end{align*}
Therefore, the total error $\delta(t)=\delta_q(t)+\delta_s(t)$ in the case when $\eps(t_k)\le 2M[k]$ is overbounded as:
\begin{align}
\label{B1}
|\delta (t)|\le M[k]+ \sup\limits_{t_k\le s\le t}|\dot{\eps}(s)|\cdot (t-t_k).
\end{align}

A difficulty for evaluation of the synchronization error is in the dependence of right-hand side of \eqref{B1}  on the trajectory of the system which is not known {\it a priori}. Therefore in order to evaluate the error dynamics we need to analyze dynamics of $\delta(t)$ and $e(t)$ simultaneously.

In order to analyze the synchronization error we make two assumptions:
\begin{enumerate}
\item[A1.] Nonlinearity $\ph(y)$ is Lipschitz continuous:
\begin{align}
\label{fr4}
|\ph(y_1)-\ph(y_2)|\le L_\ph|y_1-y_2|
\end{align}
for all $y_1$, $y_2$ and some $L_\ph>0$.
\item[A2.] The linear part of \eqref{2} is strictly passifiable:
according to the Passification Theorem \cite{Fradkov74,FradkovMN99} (see Appendix \ref{sec:pasth}), this means that the numerator $\beta(\lambda)$ of the transfer function $W(\lambda)=C(\lambda{\rm I}-A)^{-1}B$ $=\beta(\lambda)/\alpha(\lambda)$ is 
Hurwitz (stable) polynomial of degree $n-1$ with positive coefficients (the so-called hyper-minimum-phase (HMP) property).
\end{enumerate}

It follows from condition A2 and the Passification Theorem  (see Appendix \ref{sec:pasth}), that the stability degree $\eta_0$ of the polynomial $\beta(\lambda)$ (a minimum distance from its roots to the imaginary axis) is positive and for any $\eta$: $0<\eta<\eta_0$ there exist a positive definite matrix $P=P\trn>0$ and a number $K$ such that the following matrix relations hold:
\begin{align}
\label{**}
PA_K+A_K\trn P\le -2\eta P,~~PB=C\trn,~~A_K=A-BKC.
\end{align}
Any sufficiently large real number can be chosen as the value of $K$.

The main result of this Section is formulated as follows.

{\it Theorem 1.} Let A1, A2 hold, the coder parameters $T_s$, $\rho$ and the auxiliary parameters $\eta$, $\eta'$ be chosen in order to meet the inequalities 
\begin{align}
\label{t1}
T_sb_0<1,\quad q<\rho<1,\quad 0<\eta'<\eta<\eta_0,
\end{align}
where
\begin{align*}
&q=\max\bigg\{\exp(-\eta' T_s)+\dfrac{T_sa_0(K+L_\ph)\big(1-\exp(-\eta'T_s)\big)}{2\eta'(1-T_sb_0)\sqrt{\lambda_{\min}}},\notag\\
&\dfrac{T_sa_0(K+L_\ph)}{ 2(\eta-\eta')(1-T_sb_0)\sqrt{\lambda_{\min}}}\bigg\},
\end{align*}
$\eta_0$ is a stability degree of the polynomial $\beta(\lambda)$.
Let the number $K$ is chosen from \eqref{**} and $\lambda_{\min}$ is the minimum eigenvalue of the matrix $P$ from \eqref{**}.

Let the coder range $M[k]$ be chosen as follows
\begin{align}
\label{tmk}
M[k]= M_0\rho^k,
\end{align}
where 
\begin{align*}
&M_0=W_0\dfrac{(\rho-q)(1-T_sb_0)\sqrt{\lambda_{\min}}}{r(K+L_\ph)},\\
&r=\max\bigg\{\dfrac{1-\exp(-\eta'T_s)}{2 \eta'},~~\dfrac{1}{2(\eta-\eta')}\bigg\}.
\end{align*}

Then for all $e(0)$ such that $e(0)\trn P e(0)<W_0$
the current values of the synchronization errors decrease exponentially: 
\begin{align}
\label{B4}
|\eps[k]|\le \|e[k]\|\le 2\rho^k M_0.
\end{align}

{\it Proof.} The key point of the proof is comparison of the hybrid system in question with an auxiliary continuous-time system (the {\it continuous model}) possessing useful stability and passivity properties \cite{DerevitskyFradkov74,DerevitskyFradkov81}.

 Rewrite the error equation in the following form:
\begin{align}
\dot e=Ae+B\zeta(\eps,t)-Bu,\quad \eps=Ce,
\label{01}
\end{align}
where $\zeta(\eps,t)=\ph\big(y(t)\big) -\ph\big(y(t)-\eps\big)$ satisfies the inequality $|\zeta(\eps,t)|\le L_\ph|\eps|$.

Substituting \eqref{02} into \eqref{01} we obtain
\begin{align}
\dot e=A_Ke+B\zeta(\eps,t)+BK\delta (t),
\label{03}
\end{align}
where $A_K=A-BKC$. Employing the HMP condition and the Passification Theorem, see Appendix \ref{sec:pasth}, 
pick up the $(n\times n)$-matrix $P=P\trn>0$ and the positive number $K$ such that $PA_K+A_K\trn P\le -2\eta P$, $PB=C\trn$, and choose the Lyapunov function candidate $V(e)=\dfrac{1}{2}e\trn Pe$. Introducing a new nonlinearity $\xi(\eps,t)=\zeta(\eps,t)+L_\ph\eps$, satisfying the sector condition $\xi\eps\ge 0$ and making the change $K\to K+L_\ph$, transform equation \eqref{03} to the form
\begin{align}
\dot e=A_Ke+B\xi-B(K+L_\ph) \delta(t).
\label{04}
\end{align}
The time derivative of $V(e)$ is evaluated as follows:
\begin{align*}
\dot V=e\trn\big(PA_K+A_K\trn P\big)e-e\trn PB\big(K+L_\ph\big)\delta(t),
\end{align*}
or, after simple algebra
\begin{align}
\dot V\le -2\eta V+|\eps|(K+L_\ph)|\delta(t)|,
\label{05}
\end{align}
where the equivalent disturbance $\delta(t)$ satisfies the inequality
\begin{align*}
|\delta(t)|\le M[k]+ \sup\limits_{t_k\le \tau\le t_{k+1}}\|\dot Ce(\tau)\|\cdot T_s.
\end{align*}

Overbounding $\|\dot e(\tau)\|$ from \eqref{04} and taking into account that $\|e\|\le \sqrt{V(e)/\lambda_{\min}}$, where $\lambda_{\min}$ is the minimum eigenvalue of $P$ yields 
\begin{align}
|\delta(t)|\le M[k]+T_s\left(a_0\sqrt{\bar V[k+1]}+b_0\bar \delta[k+1]\right),
\label{06}
\end{align}
where $\bar V[k+1]=\sup\limits_{t_k\le t\le t_{k+1}}V\big(e(t)\big)$,
$\bar \delta[k+1]=\sup\limits_{t_k\le t\le t_{k+1}}|\delta(t)|$, $a_0=\big(\|CA_K\|$$+|CB|\cdot\|C\|\big)/\sqrt{\lambda_{\min}}$,
\par\noindent$b_0=|CB|\cdot\left(K+L_\ph\right)$.

Let $T_sb_0<1$, i.e. $T_s<\left(|CB|\cdot\left(K+L_\ph\right)\right)^{-1}$. Taking the supremum of the left part of \eqref{06} over $[t_k, t_{k+1}]$ we arrive at an inequality
\begin{align}
\bar\delta[k+1]\le \dfrac{M[k]}{1-T_sb_0}+\dfrac{T_sa_0}{1-T_sb_0}\sqrt{\bar V[k+1]}.
\label{07}
\end{align}

Substituting \eqref{07} into \eqref{05} we obtain the functional-differential inequality for $t_k\le t \le t_{k+1}$
\begin{align}
\dot V\le -2\eta V+\sqrt{V\big(e(t)\big)}\big(a\sup\limits_{t_k\le t'\le t}\sqrt{V(t')}+b[k]\big),~~k=0,1,\dots
\label{08}
\end{align}
where 
\begin{align*}
&a=\dfrac{T_sa_0(K+L_\ph)}{(1-b_0T_s)\sqrt{\lambda_{\min}}},~~b[k]= \dfrac{M[k] (K+L_\ph)} {(1-b_0T_s)\sqrt{\lambda_{\min}}}.
\end{align*}

To estimate solutions of \eqref{08} employ Lemma 1 (see Appendix \ref{proofl1}) with $V(t)=V\big(e(t)\big)$. Choosing the coder parameters satisfying conditions of the Theorem and checking the conditions of Lemma 1 complete the proof.

Based on Theorem 1 the following design method is proposed. The transmission rate $T_s$ should be chosen based on condition A1. Parameter $\rho$ in \eqref{mk} should be found independently on $T_s$ based on the accessible stability degree $\eta$ of the continuous model \eqref{04} (the maximal value of $\eta$ corresponds to the stability degree $\eta_0$ of the numerator $\beta(\lambda)$ of the transfer function $W(\lambda)=C(\lambda{\rm I}-A)^{-1}B$. Finally, the value of $\eta$ should be found based on the solution of LMI \eqref{**} for the chosen control gain $K$. 

{\it Remark 1}. Conditions $q<1$, $T_sb_0<1$ are always true for sufficiently small $T_s>0$ (i.e. for sufficiently large capacity of the communication channel). Indeed, the following relations hold up to the second terms in $T_s$:
\begin{align*}
q\approx \max\bigg\{1-\eta'T_s+\dfrac{ T_s ^2a_0(K+L_\ph)}{2\sqrt{\lambda_{\min}}},~~\dfrac{T_sa_0(K+L_\ph)}{2(\eta-\eta')\sqrt{\lambda_{\min}}}\bigg\}.
\end{align*}

The threshold for $T_s$ will be 
\begin{align}
\label{***}
&T_s<\dfrac{1}{K+L_\ph}\min      \Big\{ 2\eta'\sqrt{\lambda_{\min}}\cdot a_0^{-1},
\notag\\&2(\eta-\eta')\sqrt{\lambda_{\min}}\cdot a_0^{-1},~~|CB|^{-1}\Big\}.
\end{align}

It follows from \eqref{***} that grows of the controller gain $K$ leads to increase of the transmission rate required for decay of the synchronization error.

The dependence 
\eqref{B4} will be used for numerical analysis in Section~\ref{Sec:example}.

{\it Remark 2.} In stochastic framework the estimates of the mean square value of the synchronization error can be obtained. There is a significant body of work in which the quantization error signal $\delta (t)$ is modeled as an extra additive white noise. This assumption, typical for digital filtering theory, is reasonable if the quantizer resolution is high \cite{Curry70}, but it needs modification for the case of a low number of quantization levels \cite{NairFagnani07}.

{\it Remark 3.} For practice, it is reasonable to choose the  coder range $M[k]$ separated from zero. The following zooming strategy for a quantizer range may be recommended instead of \eqref{mk}:
\begin{align}
\label{mklim}
M[k]=(M_0- M_{\infty})\rho^k+ M_{\infty}, ~~k=0,1,\dots ,
\end{align}
where  $0<M_{\infty}<M_0$ stands for the limit value of $M[k]$.

\section{Example. Synchronization of the chaotic Chua systems}\label{Sec:example}

Let us apply the above results to synchronization of two chaotic Chua systems coupled via a channel with limited capacity.

{\it Master system}. Let the master system \eqref{2} be represented by the following {\it Chua system}:
\begin{align}
\label{chuagen}
&\begin{cases}
\dot x_1=p(-x_1+\ph(y_1)+x_2),\quad t\ge 0,\cr
\dot x_2=x_1-x_2+x_3\cr
\dot x_3=-qx_2,
\end{cases}\\
&y_1(t)=x_1(t),\nonumber
\end{align}
where $y_1(t)$ is the master system output, $p$, $q$ are 
known parameters, $x=[x_1,x_2,x_3]\trn\inr^3$ is the state vector; $\ph(y_1)$ is a piecewise-linear function,  having the form:
\begin{align}
\ph(y)&= m_0y+m_1(|y+1|-|y-1|),
\label{phi}
\end{align}
where $m_0$, $m_1$ are given parameters.

Evidently, Chua system \eqref{chuagen} may be represented in Lurie form \eqref{2} with the matrices:
\begin{align}
\label{abc}
A=\begin{bmatrix}
-p&p&0\cr
1&-1&1\cr
0&-q&0
\end{bmatrix},\quad 
B=\begin{bmatrix}
p\cr
0\cr
0
\end{bmatrix},\quad C=[1,0,0].
\end{align}
It is easy to check that the linear part of Chua system satisfies the HMP condition. Indeed, for the triple $(A,B,C)$ from \eqref{abc} the transfer function $W(\lambda)=C(\lambda{\rm I}-A)^{-1}B=\beta(\lambda)/\alpha(\lambda)$ is as follows:
\begin{align*}
W(\lambda)=\dfrac{p(\lambda^2+\lambda+q)}
{\lambda^3+(1+p)\lambda^2+q\lambda+pq}.
\end{align*}
The numerator $\beta(\lambda)= p(\lambda^2+\lambda+q)$ is a Hurwitz polynomial of degree $2$, i.e. the HMP condition holds for all $p>0$, $q>0$.

{\it Slave system}. 
Correspondingly, the slave system equations \eqref{1} for the considered case becomes
\begin{align}
\label{chuares}
&\begin{cases}
\dot z_1=p\big(-z_1+\ph(y_2)+z_2+u(t)\big),\quad t\ge 0,\cr
\dot z_2=z_1-z_2+z_3\cr
\dot z_3=-qx_2,
\end{cases}\\
&y_2(t)=z_1(t),\nonumber
\end{align}
where $y_2(t)$ is the slave system output, $z=[z_1,z_2,z_3]\trn\inr^3$ is the state vector, $\ph(y_2)$ is defined by \eqref{phi}. 

{\it Controller} has a form \eqref{3bar}, where the control gain $K$ is a design parameter.

{\it Coding procedure} has a form \eqref{mk}, \eqref{qu11}. The input signal of the coder is $\eps(t)$. The error signal $\bar \eps (t)$ of the controller \eqref{3bar} is found by holding the value of $\bar \eps [k]$ over the sampling interval $[kT_s,(k+1)T_s)$, $k=0,1,\dots$. The initial value $M_0$ of the coder range and the decay factor $\rho$ in  \eqref{mk} are design parameters.

\begin{figure}[htpb]
\centering
\includegraphics[width=85mm] {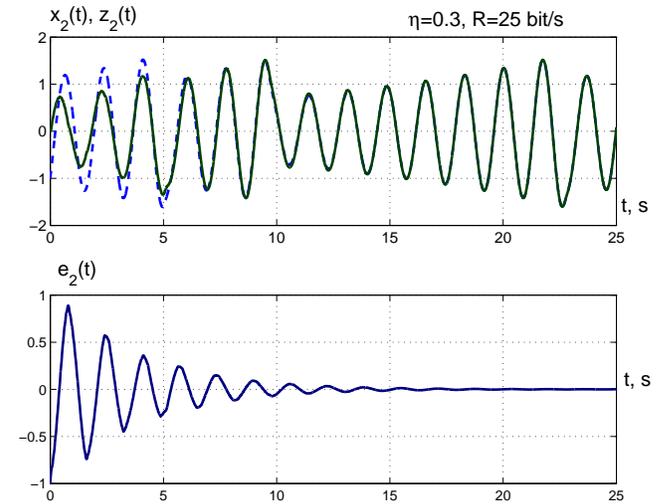}
\caption{Time histories: $x_2(t)$ (dash line), $z_2(t)$ (solid line) and synchronization error $e_2(t)=x_2(t)-z_2(t)$ for $\eta=0.3$, $ R=25 $~bit/s.}
\label{x2z2eta03r25}
\end{figure}

The following parameter values were chosen for simulation:
\begin{itemize}
\item[--] Chua system parameters: $p=10$, $q=15.6$, $m_0=0.33$, 
$m_1=0.945$; 
\item[--] the controller gain $K=10$. Feasibility of relations \eqref{A1a} for this value of $K$ and the given matrices $A$, $B$, $C$ is checked by means of {\sc Yalmip} package \cite{yalmipconf};
\item[--] the sampling time $T_s$ was taken from the interval $T_s\in[0.02,0.1]$~s for different simulation runs (a corresponding interval for the transmission rate $R$ is $R\in[10,50]$~bit/s);
\item[--] the coder parameters: $M_0=5$, $\rho=\exp(-\eta T_s)$, $\eta=0.3$;
\item[--] the initial conditions for the master and slave systems were: $x =[3,-1,0.3]\trn$, $z= 0$;
\item[--] the simulation final time $t_\text{fin}=1000$~s.
\end{itemize}
The normalized state synchronization error
\begin{align}
Q=\dfrac{\max\limits_{ 0.8t_\text{fin}\le t \le t_\text{fin}}\|e(t)\|}{\max\limits_{0\le t \le t_\text{fin}}\|x(t)\|},
\label{Q}
\end{align} 
where $\delta_y(t)= y_1(t)- \bar {y}_1 (t)$, $e(t)=x(t)-z(t)$ was calculated.  

\begin{figure}[htpb]
\centering
\includegraphics[width=85mm] {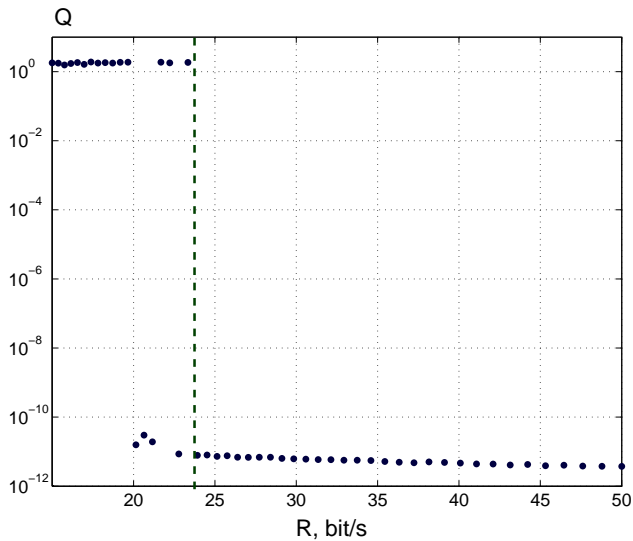}
\caption{Normalized synchronization error $Q$  v.s. transmission rate $ R $.}
\label{QotR03}
\end{figure}

Simulation results are plotted in Figs.~\ref{x2z2eta03r25},\ref{QotR03}. 

Synchronization performance may be evaluated based on time histories of the state variables $x_2(t)$, $z_2(t)$ and the synchronization error $e_2(t)$. Typical trajectories are depicted in Fig.~\ref{x2z2eta03r25}. As seen from the plots, the synchronization transient time is about $15$~seconds, which agrees with the chosen value of the coder parameter $\eta$.

The logarithmic graph of the normalized synchronization error $Q $  as a function of the transmission rate $ R $ is shown in Fig.~\ref{QotR03}. It is seen from this plot that if the transmission rate exceeds the {\it minimal bound } $R_{\min}\approx 23$~bit/s, the proposed controlled synchronization strategy ensures asymptotical vanishing the synchronization error. If the transmission rate is less that the bound $R_{\min}$, the synchronization is not always possible.

{\it Remark 4}. An idealized problem has been considered in this paper to highlight the effect of the data-rate limitations in the closed-loop synchronization of nonlinear systems. In real-world problems external disturbances, measuring errors and channel imperfections should be taken into account. Evidently in the presence of irregular nonvanishing, asymptotic convergence of the master and slave systems trajectories cannot be achieved.

{\it Remark 5}. Similar results are obtained if the control signal is also subjected to information constraints.

{\it Remark 6}. Optimality of the binary coder for synchronization under information constraints was established in \cite{FradkovAndrievskyEvans_PRE06} for the case when the master system output $y_1(t)$ rather than the output synchronization error $\eps(t)$ is transmitted over the channel. The problem of coder optimization for the considered case is under investigation.

{\it Remark 7}. The results obtained for synchronization problem can be extended to tracking problems in a straightforward manner, if the reference signal is described by an {\it external} ({\it exogenious}) state space model.
 
\section{Conclusion}
Limit possibilities of controlled synchronization systems under
information constraints imposed by limited information capacity of
the coupling channel are evaluated. It is shown that the framework proposed in \cite{FradkovAndrievskyEvans_PRE06}, is suitable not only for observer-based synchronization but also for controlled master-slave synchronization via a communication channel with limited information capacity.

 We propose a simple coder-decoder scheme and
 provide theoretical analysis for
multi-dimensional master-slave systems represented in Lurie
form. 
An output feedback control law is proposed based on the Passification Theorem \cite{Fradkov74,FradkovMN99}. It is shown that the synchronization error exponentially tends to zero for  sufficiently high
 transmission rate (channel capacity). The key point of the synchronization analysis is comparison of the hybrid system in question with an auxiliary continuous-time system ({\it the continuous model}) possessing useful stability and passivity properties. Such an approach was systematically developed in the 1970s under the name of {\it the Method of Continuous Models} \cite{DerevitskyFradkov74,DerevitskyFradkov81}.

 The results are
applied to controlled synchronization of two chaotic Chua systems
 via a communication channel with limited capacity. 
Simulation results illustrate and confirm the theoretical analysis. 

Unlike many known papers on control of nonlinear systems over a limited-band communication channel, we propos and justify a simple coder/decoder scheme, which does  not require transmission of the full system state vector over the channel. A constructive design method for controller and coder/decoder pair is proposed and estimates of the convergence rate are given. The results obtained for synchronization problem can be extended to tracking problems in a straightforward manner, if the reference signal is described by an {external} ({exogenious}) state space model.

Future research is aimed at examination of more complex system configurations, where  channel imperfections (drops, errors, delays) will be taken into account.

\begin{appendix}

\section{Passification Theorem}\label{sec:pasth}
Consider a linear system
\begin{align}
\dot{e}=Ae+B\xi(t), \quad\eps=Ce 
\label{A1}
\end{align}
with transfer function $W(\lambda)=C(\lambda I-A)^{-1}B$$=b(\lambda)/a(\lambda)$,
where $b(\lambda)$, $a(\lambda)$ are polynomials, degree of $a(\lambda)$ is $n$, degree of  $b(\lambda)$ is not greater than $n-1$. The system is called {\it hyper-minimum phase} (HMP), if $ b(\lambda)$ is  Hurwitz polynomial of degree  $n-1$ with positive coefficients. To find existence conditions for a quadratic Lyapunov function we need the following result.

{\it Passification Theorem \cite{Fradkov74,FradkovMN99}}. There exist positive-definite matrix $P=P\trn>0$ and a number $K$ such that
\begin{align}
PA_K+A_K\trn P<0,~~ PB=C\trn,~~ A_K=A-BKC  
\label{A1a}
\end{align}
if and only if   $W(\lambda)$ is HMP.

Consider a linear system with feedback
\begin{align}
\dot{e}=A_Ke+B\xi(t),~~ \eps=Ce,~~ A_K=A-BKC.
\label{A2}
\end{align}

{Remark A1}. It follows from the Passification Theorem that there exist a quadratic form $V(e)=e\trn Pe$ and a number $K$ such that time derivative $\dot V(e)$ of $V(e)$ along trajectories of  \eqref{A2} satisfies relation
\begin{align}
&\dot V(e)<0 ~~\text{for}~~ \xi\eps\ge 0, ~x\ne 0
\label{A3}
\end{align}
if and only if $W(\lambda)$ is HMP. Indeed, assume that $K$ is fixed. Relation \eqref{A3} is equivalent to existence of the matrix $P=P\trn>0$ such that $e\trn P(A_Ke+B\xi)+\xi Ce<0$ for $x\ne 0$. Since $\xi$ is arbitrary, the latter in turn, is equivalent to matrix relations  $PA_K+A_K\trn P<0$, $PB=C\trn$ and, by Passification Theorem, to HMP condition. 

{\it Remark A2}. It also follows from Passification Theorem that if HMP condition holds then $K$ satisfying \eqref{A1a} can be chosen sufficiently large. Besides, zero matrix in the right hand side of  the inequality in \eqref{A1a} can be replaced by matrix $-\mu P$ for sufficiently small $\mu>0$.

\section{Lemma 1}\label{proofl1}
{\it Lemma 1.} 
 Consider the functional-differential inequality 
\begin{align}
&\dot V(t)\le -2\eta V(t)+\sqrt{V(t)}\bigg(a\sup\limits_{t_k\le t'\le t}\sqrt{V(t')}+b[k]\bigg),\notag\\&t_k\le t<t_{k+1},~~t_k=kT_s, ~~k=0,1,\dots.
\label{*}
\end{align}

Let $0<\eta$ and $q<1$, where 
\begin{align}
q=\max\left\{\exp(-\eta' T_s)+\dfrac{a\big(1-\exp(-\eta' T_s)\big)}
{ 2\eta'},\dfrac{a}{2(\eta-\eta')}\right\}. 
\label{19b}\end{align}

Then the following statements hold:
\begin{align}
&\text{1.} \quad W[k+1]\le qW[k]+rb[k], \notag\\ &\text{where}~~ r=\max\bigg\{\dfrac{1-\exp(-\eta' T_s)}{2 \eta'},~\dfrac{1}{2(\eta-\eta')}\bigg\};
&\label{35a}\\
&\text{2.} \quad W[k]\le q^kW_0+r\sum\limits_{i=0}^{k-1}{b_iq^{k-i-1}};
\qquad \qquad &\label{35b}\\
&\text{3.  \quad If}~~q<\rho<1~~\text{and}~b[k]=\dfrac{\rho - q}{r}W[k] \notag\\ &\quad\text{then}~~W[k]~~\text{decays exponentially:}
~~W[k]=\rho^kW_0. \qquad &\notag\end{align} 

{\it Proof of Lemma 1.}

Make the change $W(t)=\sqrt{V(t)}$ and introduce the notation $W[k]=W(t_k)$, $\overline{W}(t)=\max \limits_{t_{k-1}\le t\le t } W (t)$. Then inequality \eqref{*} becomes 
\begin{align}
\dot W\le-\eta W+\dfrac{a}{2}\overline {W}(t)+\dfrac{b[k]}{2}.
\label{010}
\end{align}

Let us prove an auxiliary statement:
for any $\eta'$: $0<\eta'<\eta$ and for any $k=0,1,\dots$ the function $W(t)$ decreases and the following inequality holds:
\begin{align}
\label{p1}
&W(t)\le \max\bigg\{W(T_k)\exp\big(-\eta'(t-t_k)\big) \notag\\ &+\dfrac{aW(t_k)+b[k]}{2\eta'}\big(1-\exp(-\eta'(t-t_k))\big),
\notag\\ &\dfrac{aW(t_k)+b[k]}{2(\eta-\eta')}\bigg\}.
\end{align}

To this end consider the inequality 
\begin{align}
\label{p2}
-\eta W(t)+\dfrac{a}{2}\overline{W}(t)+\dfrac{b[k]}{2}\le \eta' W(t).
\end{align}

It holds for $t=t_k$, and, obviously holds for some interval $t_k<t<t'$. Let $t'$ be the maximum of instants for which \eqref{p2} holds. It follows from \eqref{010} that $\dot W(t)\le 0$ and $W(t)\le W(t_k)$ for $t_k\le t<t'$. Hence, $\overline{W}(t)=W(t)$ for $t_k\le t\le t'$, and \eqref{p2} becomes
\begin{align}
\label{p3}
\dfrac{a}{2}W(t_k)+\dfrac{b[k]}{2}\le (\eta-\eta')W(t)
\end{align}
and it holds iff
\begin{align}
\label{p4}
W(t)\ge \dfrac{aW(t_k)+b[k]}{2(\eta-\eta')}.
\end{align}

Therefore, as long as \eqref{p4} holds, the value of $W(t)$ decreases exponentially and
\begin{align}
\label{p5}
&W(t)\le W(t_k)\exp\big(-\eta'(t-t_k)\big) \notag\\ &+\dfrac{aW(t_k)+b[k]}{2\eta'}
\Big(1-\exp\big(-\eta'(t-t_k)\big)\Big).
\end{align}

On the other hand, if the opposite inequality $W(t)\le \dfrac{aW(t_k)+b[k]}{2(\eta-\eta')} $ holds for some $t=t'$, it holds for all $t: t'<t\le t_{k+1}$ owing to monotone decrease of $W(t)$. Therefore, \eqref{p1}
is proven.

The first statement of the lemma follows from \eqref{p1} if we set $t=t_{k+1}$,
collect terms with $W(t_k)$ and $b[k]$ and take their maximum values. The second statement is derived by iteration of \eqref{35a} and the third one is obtained by direct substitution.

\end{appendix}

\section*{Acknolegements}
This work was supported by NICTA, University of Melbourne and The Russian Foundation for Basic Research %
(projects RFBR 05-01-00869, 06-08-01386).
\bibliography{IEEEabrv,ArxFAEContr07}
\end{document}